\documentclass[10pt]{article}
\usepackage{repstyle}
\usepackage{latexsym,epsf,amssymb,subeqn}
\usepackage{graphicx}

\newcommand{\eq}{\begin{equation}}
\newcommand{\eeq}{\end{equation}}
\newcommand{\R}{\mathbb R}
\newcommand{\C}{\mathbb C}

\newcommand{\Oh}{{\cal O}}
\newcommand{\sfrac}[2]{\mbox{\large{$#1\over#2$}}}
\newcommand{\mfrac}[2]{\mbox{\Large{$#1\over#2$}}}
\newcommand{\sDelta}{{\mbox{\footnotesize{$\Delta$}}}}
\newcommand{\dt}{\sDelta t}

\newcommand{\eps}{\varepsilon}

\begin{document}

\title{On the Construction of Splitting Methods by Stabilizing Corrections
with Runge-Kutta Pairs}
\author{Willem Hundsdorfer%
\footnote{
{CWI, Science Park 123, Amsterdam, The Netherlands.}
{E-mail: willem.hundsdorfer@cwi.nl}
} }
\date{}
\maketitle

\begin{abstract} \noindent
In this technical note a general procedure is described to construct
internally consistent splitting methods for the numerical solution 
of differential equations, starting from matching pairs of explicit 
and diagonally implicit Runge-Kutta methods.  
The procedure will be applied to suitable second-order pairs, and we 
will consider methods with or without a mass conserving finishing stage. 
For these splitting methods, the linear stability properties are 
studied and numerical test results are presented.

\end{abstract}

\section{Introduction}

In this note we will discuss a class of splitting methods for solving
initial value problems for ordinary differential equations (ODEs)
\eq
\label{eq:ODE}
u'(t) = F(t,u(t)) \,, \qquad u(0) = u_0 \,,
\eeq
with given $u_0\in\R^M$, $F:\R\times\R^M\rightarrow \R^M$ and
dimension $M\ge1$.
For many practical problems there is a natural decomposition
\eq
\label{eq:Deco}
F(t,u) \,=\, F_0(t,u) + F_1(t,u) + \cdots + F_s(t,u)
\eeq
in which the separate component functions $F_j$ are more simple than
the whole $F$, and where $F_0$ is a non-stiff or mildly stiff term 
that can be treated explicitly in a time stepping method. For such problems
we will study a class of stabilizing correction splitting methods, 
where explicit predictions are followed by corrections that are implicit 
in one of the $F_j$ terms, $j=1,2,\ldots,s$.
The methods will be constructed such that all intermediate stages yield
consistent approximations to the exact solution.

\subsection{Stabilizing corrections: general procedure}

Consider a pair of Runge-Kutta methods, consisting of a diagonally implicit 
method with coefficients $a_{ik}$ ($k\le i$), and an explicit method
with coefficients $\hat{a}_{ik}$ ($k<i$), and assume these two methods 
have the same abscissae $c_i = \sum_{k\le i} a_{ik} = \sum_{k<i} \hat{a}_{ik}$.
If the methods are applied to (\ref{eq:ODE}), with known $u_n \approx u(t_n)$,
$t_n = n\dt$, the $i$-th stage of the implicit method reads
\begin{subequations}
\eq
\label{eq:stageEx}
y_i \,=\, u_n \,+\, \dt \sum_{k=1}^{i} a_{ik} F(t_n+c_k\dt, y_k) \,,
\eeq
and for the explicit method it reads
\eq
\label{eq:stageIm}
y_i \,=\, u_n \,+\, \dt \sum_{k=1}^{i-1} \hat{a}_{ik} F(t_n+c_k\dt, y_k) \,.
\eeq
\end{subequations}

We will combine these methods for problems with decomposition (\ref{eq:Deco}) 
by using the explicit formula as a predictor, followed by correction steps 
for the implicit terms. The general procedure is:
\eq
\label{eq:StageSC}
\setlength{\arraycolsep}{1mm}
\left\{
\begin{array}{rcl}
x_{i,0} &=& 
\displaystyle
u_{n} \,+\, \dt \sum_{k=1}^{i-1} \hat{a}_{ik} F(t_n+c_k\dt, y_k) \,,
\\[2mm]
x_{i,j} &=& 
\displaystyle
x_{i,j-1} \,+\, \dt \sum_{k=1}^{i-1} \big(a_{ik} - \hat{a}_{ik}\big)
F_j(t_n+c_k\dt, y_k) 
\\[3mm]
&  & 
\,+\,\, \dt\,a_{ii} F_j(t_n+c_i\dt, x_{i,j}) 
\qquad\qquad (j = 1,2,\ldots,s) \,,
\\[2mm] 
y_i &=& x_{i,s} \,.
\end{array}
\right.
\eeq
The implicit stages, where the $x_{i,j}$ are computed, mainly serve
to stabilize the process, allowing the $F_j$ terms to be stiff. Since
the two Runge-Kutta methods  
have the same abscissae $c_k$, the vectors $x_{i,0}, x_{i,1},\ldots,x_{i,s}$ 
will all be consistent approximations to $u(t_n+c_i\dt)$.
Following the terminology of \cite{Ma90}, we will call (\ref{eq:StageSC})
a stabilizing correction procedure.

The best known method of this type is obtained by combining the explicit
Euler method with the implicit trapezoidal rule or the implicit Euler 
method. This method is known as the Douglas method because of the close 
relation to ADI methods developed by J.\,Douglas\,Jr.\ and co-workers 
\cite{Do62,DoGu64} for multi-dimensional parabolic problems with 
dimension splitting, cf.\ also \cite[p.\,373]{HuVe03}.  
A class of methods with two stabilizing correction stages of the form
(\ref{eq:StageSC}) has been derived in \cite{Hu02}. In the present paper 
we will consider a more general approach, starting with Runge-Kutta
pairs of order two with two or three stages.

\subsection{Outline of the paper}

In Section~2 we will derive stabilizing correction schemes based on 
suitable pairs of Runge-Kutta methods of order two. After the stabilizing 
correction stages, a finishing stage can be appended to guarantee the
preservation of linear invariants, for example mass conservation. 
  
The stability properties of the methods are examined in Section~3 for
scalar linear test equations. It will be seen that the methods with the
appended finishing stage become unstable in general for stiff problems 
with $s\ge2$.  In Section~4 some numerical test results are presented 
for a 2D reaction-diffusion problem, where we will consider
$s=1$ (splitting of reaction and diffusion) as well as $s=2$ (with 
dimension splitting). The final Section~5
contains remarks on generalizations and conclusions.

\section{Stabilizing correction methods of order two}

\subsection{Implicit and explicit Runge-Kutta pairs of order two}

As a starting point for a stabilizing correction method, one needs a
suitable pair of implicit and explicit methods. Here we consider pairs of
second-order Runge-Kutta methods. The implicit method is taken to be
diagonally implicit and stiffly accurate, with three stages and abscissae 
$0, \kappa, 1$.
For the explicit method we take a two-stage method with the same $\kappa$ 
as abscissa.  This pair can be represented in tableau form as
\eq
\label{eq:ImExRKa}
\begin{array}{c|ccc}
0  &  0  &  &
\\
\kappa  & a_{21}  &  \theta &
\\
1  &  b_1  &  b_2  &  \theta
\\
\hline
\rule{0mm}{4mm}  &  b_1  &  b_2  &  \theta
\end{array}
\hspace{2cm}
\begin{array}{c|ccc}
0  &  &  &
\\
\kappa  &  \hat{a}_{21}  &  &
\\
1  &  \hat{b}_1  &  \hat{b}_2  &
\\
\hline
\rule{0mm}{4mm}  &  \hat{b}_1  &  \hat{b}_2  &  0
\end{array}
\quad \raisebox{-1em}{,}
\eeq
where the explicit method is written with a reducible extra stage
to make it more similar to the implicit method.  This stage will 
not be used in computations. Further it will be assumed that 
\begin{subequations}
\eq
\label{eq:order2aa}
\kappa \,=\, \hat{a}_{21} \,=\, a_{21} + \theta \,,
\eeq
so that the coefficients match the abscissae. Then the conditions for 
order two are
\eq
\label{eq:order2a}
b_1 + b_2 + \theta = 1 \,, \quad 
b_2 \kappa + \theta = \sfrac{1}{2} \,, \quad
\hat{b}_1 + \hat{b}_2 = 1 \,, \quad
\hat{b}_2 \kappa = \sfrac{1}{2} \,.
\eeq
\end{subequations}
This leaves us with two free parameters.

As an alternative we will also consider an augmented explicit method 
where the finishing stage of the implicit method is copied, giving
\eq
\label{eq:ImExRKb}
\setlength{\arraycolsep}{2mm}
\begin{array}{c|ccc}
0  &  0  & &
\\
\kappa & a_{21}  &  \theta &
\\
1  &  b_1  &  b_2  &  \theta
\\
\hline
\rule{0mm}{4mm}  &  b_1  &  b_2  &  \theta
\end{array}
\hspace{2cm}
\begin{array}{c|ccc}
0  &     & &
\\
\kappa  &  \hat{a}_{21}  &   &
\\
1  &  \hat{a}_{31}  &  \hat{a}_{32}  &
\\
\hline
\rule{0mm}{4mm}  &  b_1  &  b_2  &  \theta
\end{array}
\quad \raisebox{-1em}{.}
\eeq
Together with the matching conditions
\begin{subequations}
\eq
\label{eq:order2bb}
\kappa \,=\, \hat{a}_{21} \,=\, a_{21} + \theta \,, \qquad
\hat{a}_{31} + \hat{a}_{31} = 1 \,,
\eeq
we will impose order two, leading to the conditions
\eq
\label{eq:order2b}
b_1 + b_2 + \theta = 1 \,, \quad
b_2 \kappa + \theta = \sfrac{1}{2} \,, \quad
\hat{a}_{31} + \hat{a}_{32} = 1 \,.
\eeq
\end{subequations}
This gives three degrees of freedom in the parameters.

\subsection{The stabilizing correction methods}

In the stabilizing correction stages (\ref{eq:StageSC}) the difference
between the coefficients appear. For the second-order pair 
(\ref{eq:ImExRKa}) we have
$a_{21} - \hat{a}_{21} = - \theta$, 
$b_1 - \hat{b}_1 = {\theta}/{\kappa} - \theta$
and $b_2-\hat{b}_2 = -{\theta}/{\kappa}$.
This leads to the following stabilizing correction method:
\eq
\label{eq:SCa}
\setlength{\arraycolsep}{1mm}
\left\{
\begin{array}{rcl}
v_{0} &=& u_{n} + \kappa\,\dt\,F(t_{n}, u_{n}) \,,
\\[2mm]
v_{j} &=& v_{j-1} +
\theta \dt \big(F_j(t_{n+\kappa}, v_{j}) - F_j(t_{n}, u_{n}) \big)
 \hfill (j = 1,2,\ldots,s) \,,
\\[2mm]
w_{0} &=& u_{n} + \hat{b}_1\dt\,F(t_{n}, u_{n})
+ \hat{b}_2 \dt\,F(t_{n+\kappa},v_s)  \,,
\\[2mm]
w_{j} &=& w_{j-1} +
\theta \dt \big( F_j(t_{n+1}, w_{j}) 
- \big(1 - \sfrac{1}{\kappa}\big) F_j(t_{n}, u_{n})
- \sfrac{1}{\kappa} F_j(t_{n+\kappa}, v_{s}) \big)
\\[2mm]
 &  & \quad  
\hfill (j = 1,2,\ldots,s) \,,
\\[0mm]
u_{n+1} &=& w_{s} \,.
\end{array}
\right.
\eeq
Here the $v_i \approx u(t_{n+\kappa})$ and $w_i \approx u(t_{n+1})$,
$i=0,1,\ldots,s$, are internal vectors and the intermediate time level 
is $t_{n+\kappa} = t_{n} + \kappa\dt$.  Using $\theta$ and $\kappa$ as 
free parameters we have $\hat{b}_1 = 1 - 1/(2 \kappa)$ and 
$\hat{b}_2 = 1/(2 \kappa)$.
To distinguish this method from a variant with an extra finishing stage,
to be introduced next, we will often refer to (\ref{eq:SCa}) as 
a stabilizing correction method of {\it type-{\small A}}. 

For the implicit and explicit pair (\ref{eq:ImExRKb}), the final
stage of the methods can be appended to the stabilizing correction 
stages. This leads to the following method:
\eq
\label{eq:SCb}
\setlength{\arraycolsep}{1mm}
\left\{
\begin{array}{rcl}
v_{0} &=& u_{n} + \kappa\,\dt\,F(t_{n}, u_{n}) \,,
\\[2mm]
v_{j} &=& v_{j-1} +
\theta \dt \big(F_j(t_{n+\kappa}, v_{j}) - F_j(t_{n}, u_{n}) \big)
 \hfill (j = 1,2,\ldots,s) \,,
\\[2mm]
w_{0} &=& u_{n} + \hat{a}_{31}\dt\,F(t_{n}, u_{n})
+ \hat{a}_{32} \dt\,F(t_{n+\kappa},v_s)  \,,
\\[2mm]
w_{j} &=& w_{j-1} +
\theta \dt \big( F_j(t_{n+1}, w_{j}) - \mu_1 F_j(t_{n}, u_{n})
- \mu_2 F_j(t_{n+\kappa}, v_{s}) \big)
\\[2mm]
 &  & \quad
\hfill (j = 1,2,\ldots,s) \,,
\\[2mm]
u_{n+1} &=&  u_{n} + b_1 F(t_{n}, u_{n})
+ b_2 \dt F(t_{n+\kappa},v_s) +\theta F(t_{n}, w_{s}) \,,
\end{array}
\right.
\eeq
where
\eq
\label{eq:mu}
\mu_1 \,=\, \mfrac{1}{\theta} (\hat{a}_{31}-b_1) \,, \qquad
\mu_2 \,=\, \mfrac{1}{\theta} (\hat{a}_{32}-b_2) \,.
\eeq
Note that if  $\hat{a}_{31}=\hat{b}_1$, $\hat{a}_{32}=\hat{b}_2$, 
then the formulas for the internal vectors are the same as in 
(\ref{eq:SCa}), but for this variant (\ref{eq:SCb}) there is still 
a finishing stage with the whole function $F$.
We will refer to (\ref{eq:SCb}) as a {\it type-{\small B}} stabilizing
correction method.

It will be examined how this finishing stage influences the local
accuracy and stability of the schemes. Here we can already mention
two important properties of the above methods: {internal consistency} 
for {type-{\small A}} and {type-{\small B}} methods, and 
{mass conservation} for {type-{\small B}} methods.

{\bf Internal consistency\/}: 
The internal vectors $v_j, w_j$ are consistent approximations to 
the exact solution at time levels $t_{n+\kappa}$ and $t_{n+1}$, 
respectively. This important property guarantees that steady 
state solutions of autonomous problems are returned without any 
error. Methods with this property are called {\it well-balanced}
in shallow water applications.
Many other splitting methods, based on Lie splitting or
Strang splitting, do not share this property.

{\bf Mass conservation\/}: 
Consider a {type-{\small B}} method, and suppose the ODE system 
(\ref{eq:ODE}) is such that $h^T u(t)$ is constant for any solution, 
with a weight vector $h\in\R^m$, so the system has a linear invariant. 
This is equivalent with $h^T F(t,v) = 0$ for all $v\in\R^m$. For a 
splitting of $F$ it may happen that $h^T F_j(t,v)\neq0$ for some $j$. 
In that case, due to finishing stage with the whole function $F$, 
the {type-{\small B}} method (\ref{eq:SCb}) will still preserve 
the linear invariant, $h^T u_n = h^T u_{n-1}$, but this property 
may be lost with the {type-{\small A}} method~(\ref{eq:SCa}). 
In particular for mass conservation it can be important to maintain
linear invariants in a numerical method.

If $s=1$ the above stabilizing correction methods reduce to
implicit-explicit (IMEX) methods. This special case often occurs in
practice and it will be closely examined in this paper. 
Further it is noted that the same approach with stabilizing corrections 
could be used for Runge-Kutta methods with more stages.
This will lead, however, to classes of methods with many free coefficients, 
from which it will not be easy to chose (embarrassment of riches).
Finally we mention that the two classes of splitting methods 
(\ref{eq:SCa}) and (\ref{eq:SCb}) do not include the modified
Craig-Sneyd methods constructed by in\,'t\,Hout and Welfert \cite{HoWe09}.
These methods contain explicit stages where the coefficients for $F_0$ 
are different from those for the implicit component functions $F_j$.
Comments on this are given in the last section of this note.

\subsection{Examples}

In these notes, we will focus on classes of methods that are obtained 
by specific choices of either the explicit or the implicit methods. 
In the first two examples we will have a connection between the explicit 
methods in (\ref{eq:ImExRKa}) and (\ref{eq:ImExRKb}) by taking
\eq
\label{eq:Exas}
\hat{a}_{31}=\hat{b}_1 \,, \qquad
\hat{a}_{32}=\hat{b}_2 \,.
\eeq

\begin{Exa} \label{Exa1} \rm
One of the best known explicit two-stage Runge-Kutta method is the
explicit trapezoidal rule, also known as the modified Euler method,
with coefficients
\begin{subequations}
\label{eq:Exa1}
\eq
\label{eq:Exa1a}
\kappa \,=\, 1 \,, \qquad
\hat{b}_1 \,=\, \hat{b}_2 \,=\, \sfrac{1}{2} \,.
\eeq
From the conditions for order two, it follows that the coefficients of the
corresponding implicit method are
\eq
\label{eq:Exa1b}
a_{21} \,=\, 1-\theta \,, \quad
b_1 \,=\, \sfrac{1}{2} \,, \quad
b_2 \,=\, \sfrac{1}{2} - \theta \,,
\eeq
\end{subequations}
where we will use the diagonal coefficient $\theta$ as free parameter.
The resulting stabilizing correction method (\ref{eq:SCa}) was introduced 
in \cite{Hu02}.  In this note it will be examined whether more favourable 
methods can be found within the classes (\ref{eq:SCa}) or (\ref{eq:SCb}).

\begin{figure}[b!]
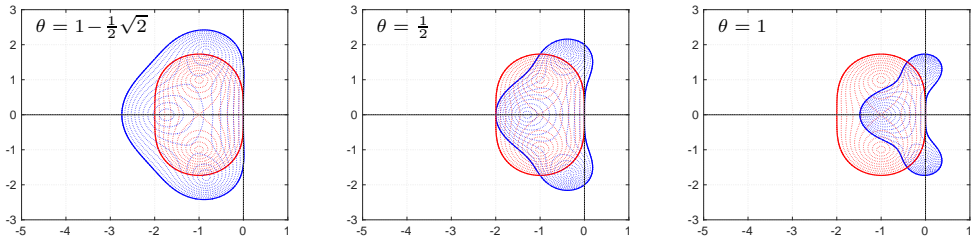

\setlength{\unitlength}{1cm}
\begin{center}
\begin{picture}(3.8,3.3)
\includegraphics[width=3.8cm]{StExa1A.eps}
\put(-3.4,2.75){\scriptsize $\theta = 1 \!-\! \frac{1}{2}\sqrt{2}$}
\end{picture}
\hspace{.5cm}
\begin{picture}(3.8,3.3)
\includegraphics[width=3.8cm]{StExa1B.eps}
\put(-3.4,2.75){\scriptsize $\theta = \frac{1}{2}$}
\end{picture}
\hspace{.5cm}
\begin{picture}(3.8,3.3)
\includegraphics[width=3.8cm]{StExa1C.eps}
\put(-3.4,2.75){\scriptsize $\theta = 1$}
\end{picture}
\caption{ \small \label{Fig:StExa1}
Stability regions of the explicit trapezoidal rule (\ref{eq:Pexa1})
indicated with red lines, and the augmented explicit method (\ref{eq:PPexa1})
with blue lines. Parameter value $\theta = 1 - \frac{1}{2}\sqrt{2}$
[left panel], $\theta = \frac{1}{2}$ [middle panel] and $\theta = 1$
[right panel].
}
\vspace{-3mm}
\end{center}
\end{figure}

To identify interesting methods, stability regions will play an important
role.  The explicit trapezoidal rule has the familiar stability function
\eq
\label{eq:Pexa1}
r_{\mbox{\scriptsize expl},A}(z) \,=\, 1 + z + \sfrac{1}{2} z^2 \,.
\eeq
For the explicit method with the extra stage in (\ref{eq:ImExRKb})
we consider (\ref{eq:Exas}), that is,
$\hat{a}_{31} = \hat{a}_{32} = \frac{1}{2}$,
which yields the stability function
\eq
\label{eq:PPexa1}
r_{\mbox{\scriptsize expl},B}(z) \,=\,
1 + z + \sfrac{1}{2} z^2 + \sfrac{1}{2}\theta z^3 \,.
\eeq
The stability function of the implicit method is given by
\eq
\label{eq:Qexa1}
r_{\mbox{\scriptsize impl}}(z) \,=\,
\frac{1 + (1-2\theta) z + (\frac{1}{2}-2\theta+\theta^2) z^2}
{(1-\theta z)^2} \;.
\eeq
This implicit method is $A$-stable if $\theta\ge\frac{1}{4}$ and it is
$L$-stable for the parameter values $\theta = 1 \pm \frac{1}{2}\sqrt{2}$.
The stability regions ${\cal E}$ of the explicit methods (\ref{eq:Pexa1}) and
(\ref{eq:PPexa1}) are presented in Figure\,\ref{Fig:StExa1} for three
parameter values.
The dotted lines are contour lines for $|r_{\mbox{\scriptsize expl},A}(z)|$
and $|r_{\mbox{\scriptsize expl},B}(z)|$ at the levels $0.1,\ldots,0.9$.
Since the implicit methods are $A$-stable, the corresponding plots
for these methods are less interesting, and therefore these are not shown.

Further we note that linearization of the diagonally implicit method 
with coefficients (\ref{eq:Exa1b}) leads to a well-known Rosenbrock-type 
method, or $W$-method \cite{HaWa96}, which
is such that its order remains two with arbitrary approximations for
the Jacobian matrix $A\approx \frac{\partial}{\partial u}F(t_n,u_n)$. 
For such methods one can apply an approximate matrix factorization
where the matrix $I - \theta\dt A$ in the Rosenbrock method is replaced
by a product $\prod_{j=1}^{s} (I - \theta\dt A_j)$ with 
$A_j\approx \frac{\partial}{\partial u}F_j(t_n,u_n)$.
The resulting method can then be viewed as a linearized version of the 
splitting method (\ref{eq:SCa}); see for instance \cite[p.400]{HuVe03}.
\hfill $\Diamond$
\end{Exa}

\begin{Exa} \label{Exa2} \rm
The choice $a_{21}=\theta$ gives a popular diagonally implicit
method where the first nontrivial stage consists of a scaled step with 
the implicit trapezoidal rule; see e.g.\ \cite[p.\,144]{HuVe03}. 
Using $\theta$ as free parameter, this implicit method has order two if
\begin{subequations}
\label{eq:Exa2}
\eq
\label{eq:Exa2a}
a_{21} \,=\, \theta \,,
\quad
\kappa = 2 \theta \,,
\quad
b_1 \,=\, \sfrac{3}{2} - \theta - \sfrac{1}{4\theta} \,,
\quad
b_2 \,=\, -\sfrac{1}{2} + \sfrac{1}{4\theta} \;.
\eeq
Then, requiring order two for the two-stage explicit method gives
the coefficients 
\eq
\label{eq:Exa2b}
\hat{b}_1 \,=\, 1 - \sfrac{1}{4\theta} \,,
\qquad
\hat{b}_2 \,=\, \sfrac{1}{4\theta} \,.
\eeq
\end{subequations}
The implicit method has again stability function (\ref{eq:Qexa1}), 
because the coefficients of $r_{\mbox{\scriptsize impl}}$ are 
determined by the order two conditions.
Likewise, the stability function of the two-stage explicit method 
(\ref{eq:ImExRKa}) is given by (\ref{eq:Pexa1}).
For the augmented three-stage explicit method we consider (\ref{eq:Exas}), 
that is, $\hat{a}_{31} = 1 - \frac{1}{4\theta}$, 
$\hat{a}_{32} = \frac{1}{4\theta}$.
It is seen by some calculations that this gives again the stability 
function (\ref{eq:PPexa1}).

So, even though the methods are different if $\theta\neq\frac{1}{2}$, 
the stability functions are the same as in the previous example and the
stability regions are as in Figure\,\ref{Fig:StExa1} for $\theta = \frac{1}{2}$
and $\theta = 1 - \frac{1}{2}\sqrt{2}$. We will see in the next section that 
this is a common property of the methods based on second-order pairs 
(\ref{eq:ImExRKa}) as well as for the related methods (\ref{eq:ImExRKb}) 
with coefficients specified by (\ref{eq:Exas}).
\hfill $\Diamond$
\end{Exa}

\begin{Exa} \label{Exa3} \rm
The implicit method of the previous example was used to construct
implicit-explicit (IMEX) methods in \cite{GKC13,GhCo16}. In these
references, the implicit method (\ref{eq:Exa2a}) was combined 
with the three-stage explicit methods (\ref{eq:ImExRKb}) with parameters
\begin{subequations}
\label{eq:Exa3}
\eq
\label{eq:Exa3a}
\theta = 1 - \sfrac{1}{2}\sqrt{2} \,, \qquad
\hat{a}_{31} = \sfrac{1}{2} - \omega \,, \qquad
\hat{a}_{32} = \sfrac{1}{2} + \omega \,.
\eeq
The values $\omega = 0$ and $\omega = \sfrac{1}{3}\sqrt{2}$ correspond to the 
choices made in \cite{GKC13} and \cite{GhCo16}, respectively. The other 
parameters are
\eq
\label{eq:Exa3b}
a_{21} = 1 - \sfrac{1}{2}\sqrt{2} \,, \qquad
\kappa = 2 - \sqrt{2} \,, \qquad
b_1 = \sfrac{1}{4}\sqrt{2} \,, \qquad
b_2 = -\sfrac{1}{4}\sqrt{2} \,,
\eeq
\end{subequations}
in agreement with (\ref{eq:Exa2a}).  Plots of the stability regions 
of these explicit method can be found in Figure\,\ref{Fig:StExa3s1}.

The IMEX methods of \cite{GKC13,GhCo16} can be obtained from (\ref{eq:SCb}) 
with $s=1$.  We will also study these methods for the case $s>1$, but it 
will be seen that stability then becomes problematic.  Finally we note 
that the type-{\small A} methods with $\hat{b}_j = \hat{a}_{3j}$, $j=1,2$, 
are less interesting in this example since the two-stage explicit methods 
are then only of order one, except for the choice 
$\omega = \frac{1}{4}\sqrt{2}$.
\hfill $\Diamond$
\end{Exa}

For the implicit method (\ref{eq:Exa2a}), used in the above examples, 
the order will be three if $\theta = \frac{1}{2} \pm \frac{1}{6}\sqrt{3}$. 
This can be interesting for IMEX applications where the dominant error 
is caused by the implicit term.

\section{Linear stability properties}

\subsection{Stability functions}

Stability and convergence will be analyzed for linear systems of
differential equations where $F_j(t,u) = A_j u + g_j(t)$.
As a first step we consider the test equation
\eq
\label{eq:testODE}
u'(t) \,=\, (\lambda_0 + \lambda_1 + \ldots + \lambda_s\big)  u(t) \,.
\eeq
Let $z_j = \dt \lambda_j$. The
stabilizing correction methods will then give a relation 
$u_{n+1} = r(z_0,z_1,\ldots,z_s) u_n$. Similar as for Runge-Kutta
methods, such a function $r$ will be called the stability function.
It will be seen that the variables $z_0,z_1,\ldots,z_s$ appear in the
stability functions only in the combinations
\eq
\label{eq:zpomega}
z = z_0 + z_1  + \ldots + z_s \,, \qquad
\varpi \,=\, \prod_{j=1}^{s} (1-\theta z_j) \,.
\eeq 

To have a clear distinction between the stabilizing correction methods
of type-{\small A} and type-{\small B}, we will use sub-indices {\small A} or 
{\small B} for these stability functions.
First we will derive a relation $w_s = q(z_0,z_1,\ldots,z_s) u_n$.
For the type-{\small A} method (\ref{eq:SCa}) this will already provide the 
stability function.

Consider the stabilizing correction methods (\ref{eq:SCa}) and
(\ref{eq:SCb}). To get the same notation we set $\mu_1 = 1 - 1/\kappa$,
$\mu_2 = 1/\kappa$ for (\ref{eq:SCa}), and $\hat{b}_1 = \hat{a}_{31}$,
$\hat{b}_2 = \hat{a}_{32}$ for (\ref{eq:SCb}).
Then, application to the test equation gives
$$
\setlength{\arraycolsep}{1mm}
\begin{array}{rcl}
v_{0} &=& u_{n} + \kappa\,z\,u_{n} \,,
\\[1mm]
v_j &=& v_{j-1} + \theta z_j (v_j - u_n) \qquad
(j=1,\ldots,s) \,,
\\[1mm]
w_0  &=& u_n + \hat{b}_1 z\,u_n + \hat{b}_2 z\,v_s \,,
\\[1mm]
w_j  &=&  w_{j-1} + \theta z_j (w_j - \mu_1 u_n - \mu_2 v_s) \qquad
(j=1,\ldots,s) \,.
\end{array} 
$$
To derive a suitable expression for the stability function, it is 
convenient to introduce $\bar{v}_j = v_j - u_n$ and 
$\bar{w}_j = w_j - \mu_1 u_n - \mu_2 v_s$. Then
$$
\begin{array}{ll}
\bar{v}_{0} \,=\, \kappa\,z\,u_{n} \,,
\qquad
&
\bar{v}_j \,=\, \sfrac{1}{1-\theta z_j} \; \bar{v}_{j-1} \,,
\\[1mm]
\bar{w}_0  \,=\, z\, u_n - (\mu_2 - \hat{b}_2 z) \bar{v}_s \,,
\qquad
&
\bar{w}_j  \,=\, \sfrac{1}{1-\theta z_j} \; \bar{w}_{j-1} \,.
\end{array}
$$
Hence $\bar{v}_s = \bar{v}_0/\varpi$ and $\bar{w}_s =  \bar{w}_0/\varpi$.
Combining these relations and using $b_1+b_2+\theta=1$ gives
$\bar{w}_s = \frac{1}{\varpi} z u_n - 
\frac{1}{\varpi^2}(\mu_2-\hat{b}_2 z) \kappa z u_n $,
which finally leads to $w_s = q(z_0,\ldots,z_s) u_n$ with
\eq
\label{eq:q}
q(z_0,z_1,\ldots,z_s)  \,=\, 1 + (1 + \mu_2\kappa) \mfrac{z}{\varpi}
-  \mu_2\kappa  \mfrac{z}{\varpi^2} + \hat{b}_2 \kappa \mfrac{z^2}{\varpi^2}
\;.
\eeq

For method (\ref{eq:SCa}) the stability function is $r_{\!A} = q$.
By use of the conditions (\ref{eq:order2a}) for order two, the 
following result is obtained:

\begin{Pro}
\label{Pro:StabfTypeA}
For all pairs (\ref{eq:ImExRKa}) with (\ref{eq:order2aa}), (\ref{eq:order2a}),
the type-{\small A} method (\ref{eq:SCa}) has stability function
\eq
\label{eq:rA}
r_A(z_0,z_1,\ldots,z_s)  \,=\, 1 \,+\, 2 \frac{z}{\varpi}
\,-\,  \frac{z}{\varpi^2} \,+\, \sfrac{1}{2}\,\frac{z^2}{\varpi^2} \;.
\eeq
\end{Pro}

This expression follows immediately from (\ref{eq:q}): if (\ref{eq:order2a}) 
then $\kappa\mu_2 = 1$ and $\hat{b}_2\kappa = \frac{1}{2}$.  
Further it should be noted that the stability function does not depend 
on the parameter $\kappa$, and the other free parameter $\theta$ only 
enters through $\varpi = \prod_{j=1}^{s}(1-\theta z_j)$.

For the type-{\small B} method (\ref{eq:SCb}), with the additional finishing 
stage, application to the test equation gives 
$u_{n+1} = (1 + b_1 z)u_n + b_2 z v_s + \theta z w_s$,
where we can use the above expressions for $v_s$ and for $w_s$ with
$\hat{a}_{32} = \hat{b}_2$. After a little calculation this leads to
\eq
r_B(z_0,z_1,\ldots,z_s)  \,=\, 1 + z
+ \big(\sfrac{1}{2} + \theta \mu_2\kappa \big) \mfrac{z^2}{\varpi}
- \theta \mu_2\kappa  \mfrac{z^2}{\varpi^2}
+ \theta  \hat{a}_{32} \kappa \mfrac{z^3}{\varpi^2} \;.
\eeq
Use of the conditions (\ref{eq:order2b}) for order two gives the
following result:

\begin{Pro} 
\label{Pro:StabfTypeB}
For all pairs (\ref{eq:ImExRKb}) with (\ref{eq:order2bb}), (\ref{eq:order2b}),
the type-{\small B} method (\ref{eq:SCb}) has stability function
\eq
\label{eq:rB}
r_B(z_0,z_1,\ldots,z_s)  \,=\, 1 \,+\, z 
\,+\, \big(\sfrac{1}{2}+\nu\big) \frac{z^2}{\varpi}
\,-\, \nu \frac{z^2}{\varpi^2} 
\,+\, \big(\sfrac{1}{2}-\theta+\nu\big)\theta \frac{z^3}{\varpi^2} \,,
\eeq
where $\nu = \theta\kappa\mu_2$.
If (\ref{eq:order2a}) and (\ref{eq:Exas}) also hold, then
\eq
\label{eq:rB'}
r_B(z_0,z_1,\ldots,z_s)  \,=\, 1 \,+\, z 
\,+\, \big(\sfrac{1}{2}+\theta\big) \frac{z^2}{\varpi}
\,-\, \theta \frac{z^2}{\varpi^2} 
\,+\, \sfrac{1}{2}\theta \frac{z^3}{\varpi^2} \;.
\eeq
\end{Pro}

Formula (\ref{eq:rB}) directly follows from 
$\kappa\hat{a}_{32} = \frac{1}{2} - \theta + \theta\kappa\mu_2$,
which is a consequence of (\ref{eq:order2b}). 
Moreover, if (\ref{eq:order2a}) with $\hat{b}_2=\hat{a}_{32}$ is 
also valid, then $\nu=\theta$, giving (\ref{eq:rB'}).

Note that in these stability functions there are now terms 
$z^{k+1}/\varpi^k$, with power in the numerator higher than in the
denominator, so it not very surprising that stability is often 
harder to achieve for these type-{\small B} methods. This will be 
discussed in detail in the next sections.

\subsection{Stability domains}

In the following it will be assumed that the implicit arguments $z_j$, 
$j\ge1$, are in a wedge ${\cal W}_\alpha$ in the left-half plane,
$$
{\cal W}_\alpha \,=\, \{ \zeta \in \C : \arg(-\zeta) \le \alpha
\;\; \mbox{or} \;\; \zeta = 0 \}
$$
with angle $\alpha \in [0,\frac{1}{2}\pi]$.  So ${\cal W}_0 = \R^{-}$ 
is the non-positive real axis and ${\cal W}_{\pi/2} = \C^{-}$ is the 
left-half plane.
For a given stability function $r$ and angle $\alpha$, we will then study 
the following domains for the explicit argument $z_0$:
\eq
\label{eq:Dalpha}
{\cal D}_\alpha \,=\, \{ z_0 \in \C : \; |r(z_0,z_1,\ldots,z_s)| \le 1
\;\; \mbox{for all $z_j \in {\cal W}_\alpha$, $j=1,2,\ldots,s$} \}\,.
\eeq

\subsubsection{Necessary conditions for stability with $\alpha = 0$}

Necessary stability conditions can be obtained by studying the limit case
$z_s\rightarrow -\infty$, with the other implicit arguments
$z_1,\ldots,z_{s-1}$ real and non-positive.  
For a method with stability function $r$ it will be required that
\eq
\label{eq:StabLim}
\lim_{z_s\rightarrow-\infty}|r(z_0,z_1,\ldots,z_s)|  \,\le\, 1
\quad\mbox{for all $z_1,\ldots,z_{s-1} \in\R^{-}$} \,.
\eeq
Further we will use the notation
\eq
\label{eq:sigmarho}
\sigma_s \,=\, \sum_{j=0}^{s-1} z_j \,, \qquad
\pi_s \,=\, \prod_{j=1}^{s-1} (1 - \theta z_j)^{-1} \,.
\eeq
With real $z_1,\ldots,z_{s-1} \le 0$ we have $\pi_s\in[0,1]$,
and if $s=1$ then $\pi_s = 1$.
  
\begin{Pro} \label{Pro:TypeAlim} 
For the type-{\small A} methods (\ref{eq:SCa}), with stability function 
$r = r_A$  given by (\ref{eq:rA}), the stability condition (\ref{eq:StabLim}) 
holds iff $\theta \ge \frac{1}{4}$.
\end{Pro}

\noindent
{\bf Proof.}
For large $z_s$ we have ${z}/{\varpi} = -\pi_s/\theta + \Oh(z_s^{-1})$,
${z}/{\varpi}^2 = \Oh(z_s^{-1})$. Consequently
\eq
\label{eq:qA}
\lim_{z_s\rightarrow -\infty}r_A(z_0,\ldots,z_s) \,=\,
\phi_A(\pi_s) \,=\, 1 \,-\, \sfrac{2}{\theta} \pi_s 
\,+\, \sfrac{1}{2\theta^2} \pi_s^2 \,.
\eeq
For this limit function $\phi_A$ it is easily seen that
$$
\begin{array}{c}
\phi_A(0) = 1 \,, \qquad
\phi_A(1) = \frac{1}{\theta^2} \big(\frac{1}{2}-2\theta+\theta^2\big)  \,,
\\[2mm]
\phi'_A(\pi_s) = 0 \;\;\mbox{for}\;\; \pi_s = 2\theta \,, 
\;\;\mbox{and}\quad
\phi_A(2\theta) = -1 \,.
\end{array}
$$
Furthermore $|\frac{1}{2}-2\theta+\theta^2| \le \theta^2$ iff
$\theta \ge \frac{1}{4}$, which provides the proof.
[To be done more carefully and more clearly, separating the cases
$s=1$ and $s\ge2$.]
\hfill $\Box$

\bigskip
Observe that $\theta \ge \frac{1}{4}$ gives exactly the parameter range 
for which the implicit method is A-stable. So this is the only requirement
to be fulfilled for stability in the limit case $z_s\rightarrow -\infty$.

As can be expected from the form of the stability functions for
type-{\small B} methods, stability is more delicate for such methods. 
To formulate the result, we define
\eq
\label{eq:NesCondB2}
\phi_B(z_0) \,=\,
\sfrac{1}{\theta^2} \big(\sfrac{1}{2} - 2\theta + \theta^2\big)
+ \sfrac{1}{\theta} \big(\sfrac{1}{2} - 2\theta + \nu\big) z_0  \,.
\eeq 

\begin{Pro} \label{Pro:TypeBlim}
Consider the type-{\small B} methods (\ref{eq:SCb}) with stability 
function $r = r_B$ given by (\ref{eq:rB}).
\\
\textbf{\textit{(a)}}\,\,
If $s=1$, then the stability condition (\ref{eq:StabLim}) holds iff
$|\phi_B(z_0)| \le 1$.
\\
\textbf{\textit{(b)}}\,\,
If $s\ge2$, then the stability condition (\ref{eq:StabLim}) cannot hold. 
\end{Pro}

\noindent
{\bf Proof.}
For large $|z_s|$ we have
$$
\begin{array}{c}
z \,=\, z_s + \sigma_s \,=\, z_s \big(1 + \Oh(z_s^{-1}) \big) \,,
\\[2mm]
\mfrac{z}{\varpi} \,=\, -\sfrac{1}{\theta}\,\pi_s
\,-\, \sfrac{1}{\theta^2} (1 + \theta \sigma_s) \pi_s z_s^{-1}
\,+\,  \Oh(z_s^{-2}) \,,
\end{array}
$$
and therefore
$$
\begin{array}{c}
\mfrac{z^2}{\varpi} \,=\, -\sfrac{1}{\theta}\,\pi_s z_s
- \sfrac{1}{\theta^2} (1 + 2 \theta \sigma_s) \pi_s  + \Oh(z_s^{-1}) \,,
\\[2mm]
\mfrac{z^3}{\varpi^2} \,=\, \sfrac{1}{\theta^2} \pi_s^2 z_s 
+ \Big(\sfrac{1}{\theta^2} \sigma_s \pi_s^2 + \sfrac{2}{\theta^3}
(1 + \theta \sigma_s)\Big) + \Oh(z_s^{-1}) \,.
\end{array}
$$
Inserting these expansions in (\ref{eq:rB})  we obtain
$$
\begin{array}{l}
r_B(z_0,\ldots,z_s) \,=\,
\Big( 1 - \sfrac{1}{\theta}\big(\sfrac{1}{2}+\nu\big) \pi_s
+ \sfrac{1}{\theta}\big(\sfrac{1}{2} - \theta + \nu\big) \pi_s^2 \Big) z_s 
\\[2mm]
\qquad
+\, \Big( 1 + \sigma_s - \sfrac{1}{\theta^2}\big(\sfrac{1}{2}+\nu\big)
(1 + 2\theta\sigma_s)\pi_s - \sfrac{\nu}{\theta^2}\pi_s^2
\\[2mm]
\qquad\qquad
+\, \big(\sfrac{1}{2} - \theta + \nu\big) 
\big(\sfrac{1}{\theta}\sigma_s\pi_s^2
+ \sfrac{2}{\theta^2}(1 + \theta\sigma_s)\big) \Big) +  \Oh(z_s^{-1}) \,.
\end{array}
$$ 

If $z_s\rightarrow -\infty$ then $|r_B|$ will tend to a finite
limit value iff the first term on the right vanishes, that is,
\eq
\label{eq:NesCondB3}
1 - \sfrac{1}{\theta} \big(\sfrac{1}{2}+\nu\big) \pi_s
+ \sfrac{1}{\theta} \big(\sfrac{1}{2} - \theta + \nu\big) \pi_s^2 \,=\, 0 \,.
\eeq
If $s\ge2$, then $\pi_s$ may take on any value between~$0$ and~$1$, 
in which case this equality cannot be satisfied. On the other hand, 
if $s=1$ we simply have $\pi_s = 1$, in which case (\ref{eq:NesCondB3}) 
holds trivially, and it then also follows that
$\lim_{z_1\rightarrow-\infty} r_B(z_0,z_1) = \phi_B(z_0)$.
\hfill $\Box$

\bigskip
To establish the connection between part (a) of this proposition and 
Proposition~\ref{Pro:TypeAlim}, note that $|\phi_B(0)| \le 1$ iff 
$\theta \ge \frac{1}{4}$, which is the parameter range for which 
the implicit method is $A$-stable.
Furthermore, it is clear from the negative result in part (b) 
that the type-{\small B} are not suited for problems with $s\ge2$.

\subsubsection{Stability domains ${\cal D}_\alpha$}

We now consider the stability domains ${\cal D}_\alpha$ for the explicit
argument $z_0$.  Useful analytic results can be very hard to derive. 
The main objective of this section is the presentation and discussion 
of plots of these domains, obtained by taking for each $z_0$ a large 
number of points $z_1,\ldots,z_s$ on the boundary of the wedge 
${\cal W}_\alpha$ with a given angle $\alpha$.
The plots will mostly be presented only for the angles $\alpha = 0$ and 
$\alpha = \frac{1}{2}\pi$, with comments on the stability domains for 
the intermediate angle $\alpha = \frac{1}{4}\pi$ given in the text.

For a given stability function $r$, it will be convenient in the 
discussion to refer to the function
\eq
\label{eq:psi}
\psi_\alpha(z_0) \,=\, 
\sup_{z_1,\ldots,z_s \in {\cal W}_\alpha} |r(z_0,z_1,\ldots,z_s)| \,.
\eeq
The set ${\cal D}_\alpha$ then consists of those $z_0\in\C$ for which
$\psi_\alpha(z_0) \le 1$. In the plots of the stability domains, also
contour lines $\psi_\alpha(z_0) = c$ will be drawn, with dotted lines, for 
the contour levels $c = 0.1,0.2,\ldots,0.9$.

\begin{figure}[b!]
\setlength{\unitlength}{1cm}
\begin{center}
\begin{picture}(3.8,3.3)
\includegraphics[width=3.8cm]{StExa1s1_11.eps}
\put(-3.50,2.70){\scriptsize $\theta = 1 \!-\! \frac{1}{2}\sqrt{2}$}
\put(-3.50,2.40){\scriptsize $\alpha = 0$}
\put(-3.50,2.10){\scriptsize $s = 1$}
\end{picture}
\hspace{.5cm}
\begin{picture}(3.8,3.3)
\includegraphics[width=3.8cm]{StExa1s1_12.eps}
\put(-3.50,2.70){\scriptsize $\theta = \frac{1}{2}$}
\put(-3.50,2.40){\scriptsize $\alpha = 0$}
\put(-3.50,2.10){\scriptsize $s = 1$}
\end{picture}
\hspace{.5cm}
\begin{picture}(3.8,3.3)
\includegraphics[width=3.8cm]{StExa1s1_13.eps}
\put(-3.50,2.70){\scriptsize $\theta = 1$}
\put(-3.50,2.40){\scriptsize $\alpha = 0$}
\put(-3.50,2.10){\scriptsize $s = 1$}
\end{picture}
\\
\begin{picture}(3.8,3.3)
\includegraphics[width=3.8cm]{StExa1s1_21.eps}
\put(-3.50,2.70){\scriptsize $\theta = 1 \!-\! \frac{1}{2}\sqrt{2}$}
\put(-3.50,2.40){\scriptsize $\alpha = \frac{1}{2}\pi$}
\put(-3.50,2.10){\scriptsize $s = 1$}
\end{picture}
\hspace{.5cm}
\begin{picture}(3.8,3.3)
\includegraphics[width=3.8cm]{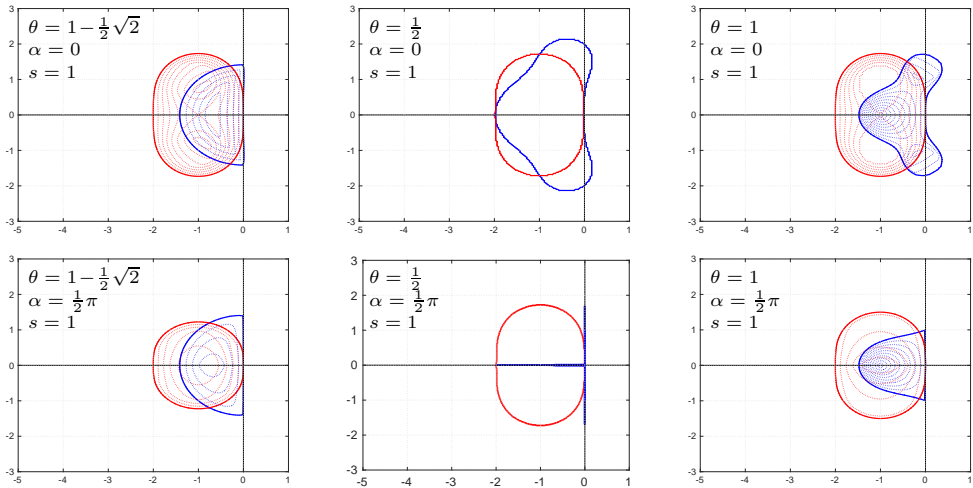}
\put(-3.50,2.70){\scriptsize $\theta = \frac{1}{2}$}
\put(-3.50,2.40){\scriptsize $\alpha = \frac{1}{2}\pi$}
\put(-3.50,2.10){\scriptsize $s = 1$}
\end{picture}
\hspace{.5cm}
\begin{picture}(3.8,3.3)
\includegraphics[width=3.8cm]{StExa1s1_23.eps}
\put(-3.50,2.70){\scriptsize $\theta = 1$}
\put(-3.50,2.40){\scriptsize $\alpha = \frac{1}{2}\pi$}
\put(-3.50,2.10){\scriptsize $s = 1$}
\end{picture}
\caption{ \small \label{Fig:StExa1s1}
Stability domains for $s=1$. Methods from Example~\ref{Exa1},~\ref{Exa2}.
Red lines for domains with $r = r_A$ [eq.~(\ref{eq:rA})],
blue for $r = r_B$ [eq.~(\ref{eq:rB'})].
From left to right: $\theta = 1 - \frac{1}{2}\sqrt{2}$
[left], $\theta = \frac{1}{2}$ [middle] and $\theta = 1$ [right].
Top row for angle $\alpha = 0$, bottom row for $\alpha = \frac{1}{2}\pi$.
}
\vspace{-3mm}
\end{center}
\end{figure}

\bigskip\noindent
{\it Methods (\ref{eq:Exa1}) with $s=1$\/}:
First we consider the methods from Example~\ref{Exa1} and
Example~\ref{Exa2}, with three values of the parameter $\theta$. 
The stability functions are $r=r_A$ and $r=r_B$, as given by the 
equations (\ref{eq:rA}), (\ref{eq:rB'}) for the type-{\small A} 
and type-{\small B} methods, respectively.  
The domains ${\cal D}_\alpha$ with angles $\alpha = 0$ and 
$\alpha = \frac{1}{2}\pi$ are shown in Figure\,\ref{Fig:StExa1s1} 
for the case $s=1$.

These domains ${\cal D}_0$ and ${\cal D}_{\pi/2}$ 
can now be compared with the stability regions ${\cal E}$
of the explicit methods, as given in Figure\,\ref{Fig:StExa1}.
It is seen that the domains ${\cal D}_0$ are equal to ${\cal E}$ for
the type-{\small A} methods. For the type-{\small B} methods this also
holds if $\theta = \frac{1}{2}$ and $\theta = 1$, but for the smallest
parameter value, $\theta = 1 - \frac{1}{2}\sqrt{2}$, the domain ${\cal D}_0$
is considerably smaller than ${\cal E}$.

The domains ${\cal D}_{\pi/2}$ are in general smaller than ${\cal E}$.
In particular, for the type-{\small B} method with parameter
$\theta = \frac{1}{2}$ the domain is reduced to a small set containing
the segment $[-2,0]$ of the negative real axis and a part of the 
imaginary axis, roughly $[-2 i, 2 i]$ (not well visible on the scale used
in these plots).
On the other hand, for corresponding type-{\small A} method with 
$\theta = \frac{1}{2}$, we again obtain the full stability region
${\cal E}$. This somewhat surprising result is quite easy to derive 
analytically, see \cite[p.\,402]{HuVe03}.

Similar plots have been made for $\alpha = \frac{1}{4}\pi$, showing 
that the stability domains for this angle are only slightly smaller 
than for $\alpha = 0$.

Finally it should be mentioned that the (dotted) contour lines 
$\psi_\alpha(z_0)=c$ with levels $c=0.1,\ldots,0.9$ are absent 
in the figure for the parameter $\theta = \frac{1}{2}$.  With this 
parameter value we always have $\psi_\alpha(z_0) \ge 1$, due to the 
fact that $|r(z_0,z_1)| \rightarrow 1$ as $z_1\rightarrow-\infty$. 
Furthermore, as a consequence, $\psi_\alpha(z_0)$ will 
not be differentiable at the boundary of the domain ${\cal D}_\alpha$, 
and this non-smooth behaviour causes the (Matlab) plotting routine
to draw staircase-shaped lines instead of smooth curves for the 
boundaries of the stability domains.

\begin{figure}[b!]
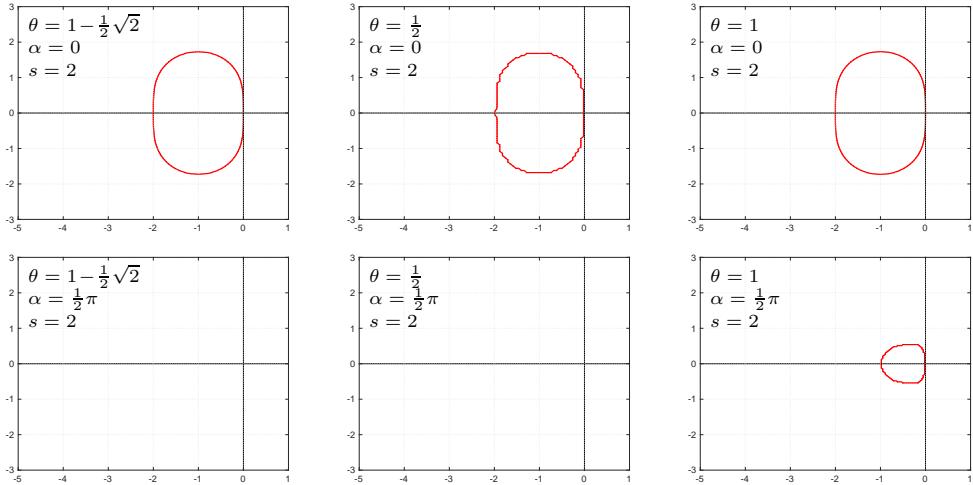

\setlength{\unitlength}{1cm}
\begin{center}
\begin{picture}(3.8,3.3)
\includegraphics[width=3.8cm]{StExa1s2_11.eps}
\put(-3.50,2.70){\scriptsize $\theta = 1 \!-\! \frac{1}{2}\sqrt{2}$}
\put(-3.50,2.40){\scriptsize $\alpha = 0$}
\put(-3.50,2.10){\scriptsize $s = 2$}
\end{picture}
\hspace{.5cm}
\begin{picture}(3.8,3.3)
\includegraphics[width=3.8cm]{StExa1s2_12.eps}
\put(-3.50,2.70){\scriptsize $\theta = \frac{1}{2}$}
\put(-3.50,2.40){\scriptsize $\alpha = 0$}
\put(-3.50,2.10){\scriptsize $s = 2$}
\end{picture}
\hspace{.5cm}
\begin{picture}(3.8,3.3)
\includegraphics[width=3.8cm]{StExa1s2_13.eps}
\put(-3.50,2.70){\scriptsize $\theta = 1$}
\put(-3.50,2.40){\scriptsize $\alpha = 0$}
\put(-3.50,2.10){\scriptsize $s = 2$}
\end{picture}
\\
\begin{picture}(3.8,3.3)
\includegraphics[width=3.8cm]{StExa1s2_21.eps}
\put(-3.50,2.70){\scriptsize $\theta = 1 \!-\! \frac{1}{2}\sqrt{2}$}
\put(-3.50,2.40){\scriptsize $\alpha = \frac{1}{2}\pi$}
\put(-3.50,2.10){\scriptsize $s = 2$}
\end{picture}
\hspace{.5cm}
\begin{picture}(3.8,3.3)
\includegraphics[width=3.8cm]{StExa1s2_22.eps}
\put(-3.50,2.70){\scriptsize $\theta = \frac{1}{2}$}
\put(-3.50,2.40){\scriptsize $\alpha = \frac{1}{2}\pi$}
\put(-3.50,2.10){\scriptsize $s = 2$}
\end{picture}
\hspace{.5cm}
\begin{picture}(3.8,3.3)
\includegraphics[width=3.8cm]{StExa1s2_23.eps}
\put(-3.50,2.70){\scriptsize $\theta = 1$}
\put(-3.50,2.40){\scriptsize $\alpha = \frac{1}{2}\pi$}
\put(-3.50,2.10){\scriptsize $s = 2$}
\end{picture}
\caption{ \small \label{Fig:StExa1s2}
Stability domains for $s=2$. Methods from Example~\ref{Exa1},~\ref{Exa2}.
Red lines for domains with $r = r_A$ [eq.~(\ref{eq:rA})]; the stability
domains for $r = r_B$ [eq.~(\ref{eq:rB'})] are empty.
From left to right: $\theta = 1 - \frac{1}{2}\sqrt{2}$
[left], $\theta = \frac{1}{2}$ [middle] and $\theta = 1$ [right].
Top row for angle $\alpha = 0$, bottom row for $\alpha = \frac{1}{2}\pi$.
}
\vspace{-3mm}
\end{center}
\end{figure}

\bigskip\noindent
{\it Methods (\ref{eq:Exa1}) with $s=2$\/}:
Next we consider the methods from Example~\ref{Exa1} and \ref{Exa2} for
the case $s=2$. For this case the stability domains are empty for the
type-{\small B} methods, in agreement with Proposition~\ref{Pro:TypeBlim}.
The domains ${\cal D}_\alpha$ for the type-{\small A} methods are presented 
in Figure\,\ref{Fig:StExa1s2} for angles $\alpha = 0, \frac{1}{2}\pi$.

With angle $\alpha = 0$ the domains still coincide with the stability
region of the underlying explicit method. 
However, for $\alpha = \frac{1}{2}\pi$ only a very small stability 
domain remains for $\theta = 1$, and for the parameters 
$\theta = 1 - \frac{1}{2}\sqrt{2}$ and $\theta = \frac{1}{2}$
the domains are empty; cf.\ \cite{LBV01} where it
was shown that ${\cal D}_{\pi/2}$ is non-empty for $s=2$ iff
$\theta \ge \frac{1}{2}+\frac{1}{6}\sqrt{3}$.

In the figure the (dotted) contour lines of $\psi_\alpha(z_0)=c$ with
$c=0.1,\ldots,0.9$ are absent. This is due to the fact that $|r(z_0,z_1,z_2)|$
tends to~$1$ as $z_1,z_2\rightarrow-\infty$, which implies that
$\psi_\alpha(z_0)\ge 1$.

The domains ${\cal D}_\alpha$ were also examined for $\alpha = \frac{1}{4}\pi$. 
For that angle, the domains with parameter values $\theta = \frac{1}{2}, 1$ 
are much larger than for $\alpha = \frac{1}{2}\pi$, namely equal to 
the stability region $\cal E$ of the corresponding explicit method, 
but for the smallest parameter $\theta = 1 - \frac{1}{2}\sqrt{2}$ the 
domain is still empty for this intermediate angle.
In fact, by considering $z_1 = \frac{2\theta-1}{2\theta^2} + i\eps$ and
$z_2\rightarrow-\infty$ it can be shown that $0\in{\cal D}_\alpha$ for
some $\alpha>0$ requires $\theta\ge\frac{1}{2}$.

\begin{figure}[h!]
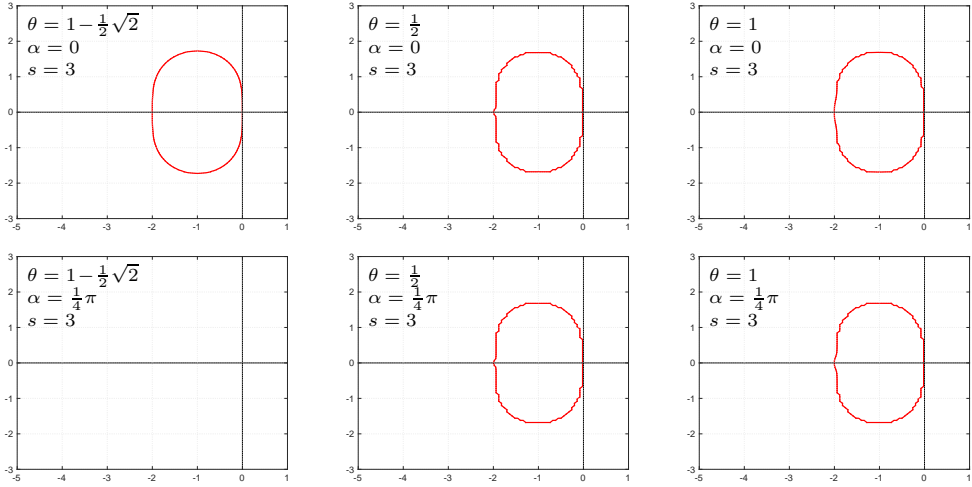

\setlength{\unitlength}{1cm}
\begin{center}
\begin{picture}(3.8,3.3)
\includegraphics[width=3.8cm]{StExa1s3_11.eps}
\put(-3.50,2.70){\scriptsize $\theta = 1 \!-\! \frac{1}{2}\sqrt{2}$}
\put(-3.50,2.40){\scriptsize $\alpha = 0$}
\put(-3.50,2.10){\scriptsize $s = 3$}
\end{picture}
\hspace{.5cm}
\begin{picture}(3.8,3.3)
\includegraphics[width=3.8cm]{StExa1s3_12.eps}
\put(-3.50,2.70){\scriptsize $\theta = \frac{1}{2}$}
\put(-3.50,2.40){\scriptsize $\alpha = 0$}
\put(-3.50,2.10){\scriptsize $s = 3$}
\end{picture}
\hspace{.5cm}
\begin{picture}(3.8,3.3)
\includegraphics[width=3.8cm]{StExa1s3_13.eps}
\put(-3.50,2.70){\scriptsize $\theta = 1$}
\put(-3.50,2.40){\scriptsize $\alpha = 0$}
\put(-3.50,2.10){\scriptsize $s = 3$}
\end{picture}
\\
\begin{picture}(3.8,3.3)
\includegraphics[width=3.8cm]{StExa1s3_21.eps}
\put(-3.50,2.70){\scriptsize $\theta = 1 \!-\! \frac{1}{2}\sqrt{2}$}
\put(-3.50,2.40){\scriptsize $\alpha = \frac{1}{4}\pi$}
\put(-3.50,2.10){\scriptsize $s = 3$}
\end{picture}
\hspace{.5cm}
\begin{picture}(3.8,3.3)
\includegraphics[width=3.8cm]{StExa1s3_22.eps}
\put(-3.50,2.70){\scriptsize $\theta = \frac{1}{2}$}
\put(-3.50,2.40){\scriptsize $\alpha = \frac{1}{4}\pi$}
\put(-3.50,2.10){\scriptsize $s = 3$}
\end{picture}
\hspace{.5cm}
\begin{picture}(3.8,3.3)
\includegraphics[width=3.8cm]{StExa1s3_23.eps}
\put(-3.50,2.70){\scriptsize $\theta = 1$}
\put(-3.50,2.40){\scriptsize $\alpha = \frac{1}{4}\pi$}
\put(-3.50,2.10){\scriptsize $s = 3$}
\end{picture}
\caption{ \small \label{Fig:StExa1s3}
Stability domains for $s=3$. Methods from Example~\ref{Exa1},~\ref{Exa2},
with the function $r = r_A$ [eq.~(\ref{eq:rA})].
From left to right: $\theta = 1 - \frac{1}{2}\sqrt{2}$
[left], $\theta = \frac{1}{2}$ [middle] and $\theta = 1$ [right].
Top row for angle $\alpha = 0$, bottom row for $\alpha = \frac{1}{4}\pi$.
}
\vspace{-3mm}
\end{center}
\end{figure}

\bigskip\noindent
{\it Methods (\ref{eq:Exa1}) with $s=3$\/}:
For general stabilizing correction methods and $s\ge3$ it can be shown that 
stability domains ${\cal D}_\alpha$ can only be non-empty if 
$\alpha \le \frac{1}{s-1}\frac{\pi}{2}$, see \cite[p.\,224]{Hu02}.
Therefore, for the case $s=3$, the plots of the stability domains 
are presented in Figure~\ref{Fig:StExa1s3} with angles $\alpha = 0$ 
and~$\frac{1}{4}\pi$.
Moreover, since we know that the type-{\small B} methods are not stable
for $s>1$, we only consider the type-{\small A} method with $r=r_A$
given by equation (\ref{eq:rA}). 

If $s=3$ then ${\cal D}_0$ equals the explicit stability region
${\cal E} = \{z_0\in\C: |1 + z_0 + \frac{1}{2} z_0|\le1\}$ for all
three parameters $\theta = 1 - \frac{1}{2}\sqrt{2}, \frac{1}{2}, 1$.
Further it is seen that ${\cal D}_{\pi/4}$ is again equal to ${\cal E}$
for $\theta = \frac{1}{2}$, it is slightly smaller 
(not well visible on this scale) for $\theta = 1$,
but for the parameter value $\theta = 1 - \frac{1}{2}\sqrt{2}$ the 
domain is now empty.
These experimental findings are in agreement with those of 
\cite[Fig.\,1,2]{Hu02}.

\bigskip\noindent
{\it Methods (\ref{eq:Exa1}) with $s>3$\/}: 
For the type-{\small A} methods with $s>3$ the above mentioned angle
bound $\alpha \le \frac{1}{s-1}\frac{\pi}{2}$ restricts the eigenvalues
of the (linearized) implicit terms. We do have the following result
for arbitrary $s$ with $\alpha = 0$:

\newpage
\begin{Pro}
Consider a type-{\small A} method with stability function $r = r_A$  
given by (\ref{eq:rA}), and suppose $\theta \ge\frac{1}{4}$, $s \ge 1$ 
and $z_0 = 0$. Then
\eq
\label{eq:0inD0}
|r(0,z_1,z_2,\ldots,z_s)| \,\le\, 1
\quad\mbox{for all $z_1,\ldots,z_{s} \in\R^{-}$} \,.
\eeq
\end{Pro}

\noindent
{\bf Proof.}
Since $z_0=0$ and the other $z_j\le0$, we have $z\le0$ and 
$\varpi \ge 1 - \theta z \ge 1$.  
For the stability function (\ref{eq:rA}), with dependence on 
$z_0,\ldots,z_s$ through $z$ and $\varpi$, we can write
$$
r \,=\, 1 \,+\, \mfrac{z}{\varpi^2} \big(2 \varpi -1\big)
\,+\, \sfrac{1}{2}\,\mfrac{z^2}{\varpi^2} \,.  
$$
Therefore $r\le1$ iff
$$
a \le 0 \qquad\mbox{with}\qquad
a = z (2\varpi - 1) + \sfrac{1}{2} z^2 \,.
$$
Since $a \le z\big(2(1-\theta z) - 1\big) + \frac{1}{2} z^2
= z + (\frac{1}{2} - 2\theta) z^2$ it follows that $a\le0$ whenever
$\theta\ge \frac{1}{4}$.  
Further we have $r\ge-1$ iff
$$
b \ge 0 \qquad\mbox{with}\qquad
b =  2\varpi^2 + z \big(2 \varpi -1\big) 
+ \sfrac{1}{2}\,z^2 \,.
$$
Since $b = \frac{1}{2}(2\varpi + z)^2 - z$, it is seen that $b\ge0$
irrespective of the value of $\theta$.
\hfill $\Box$

\bigskip
The above result is equivalent to the statement that $0 \in {\cal D}_0$ for 
any $s\ge1$.
Such a result for $z_0=0$ will provide stability on finite intervals for
non-stiff explicit terms, for which we have $|z_0| \le \dt\,L$ with a fixed 
Lipschitz constant $L$.

\begin{figure}[b!]
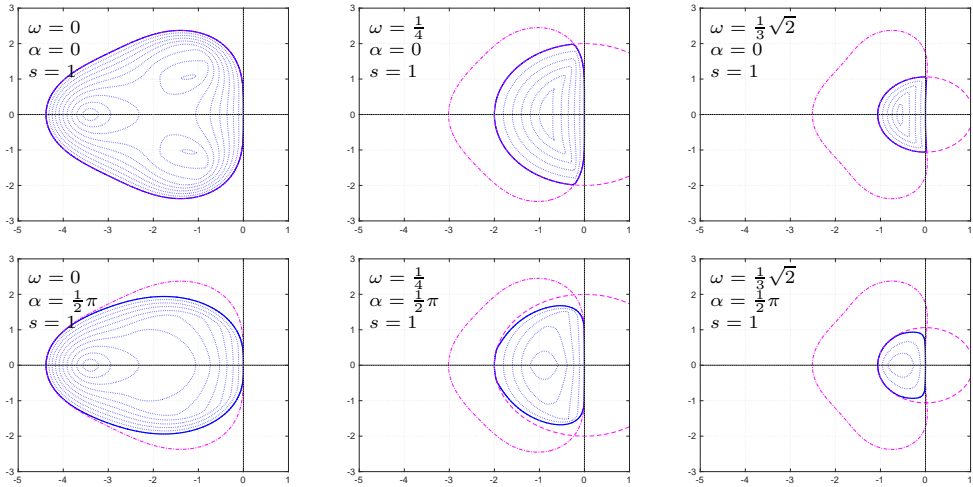

\setlength{\unitlength}{1cm}
\begin{center}
\begin{picture}(3.8,3.3)
\includegraphics[width=3.8cm]{StExa3s1_11.eps}
\put(-3.50,2.70){\scriptsize $\omega = 0$}
\put(-3.50,2.40){\scriptsize $\alpha = 0$}
\put(-3.50,2.10){\scriptsize $s = 1$}
\end{picture}
\hspace{.5cm}
\begin{picture}(3.8,3.3)
\includegraphics[width=3.8cm]{StExa3s1_12.eps}
\put(-3.50,2.70){\scriptsize $\omega = \frac{1}{4}$}
\put(-3.50,2.40){\scriptsize $\alpha = 0$}
\put(-3.50,2.10){\scriptsize $s = 1$}
\end{picture}
\hspace{.5cm}
\begin{picture}(3.8,3.3)
\includegraphics[width=3.8cm]{StExa3s1_13.eps}
\put(-3.50,2.70){\scriptsize $\omega = \frac{1}{3}\sqrt{2}$}
\put(-3.50,2.40){\scriptsize $\alpha = 0$}
\put(-3.50,2.10){\scriptsize $s = 1$}
\end{picture}
\\
\begin{picture}(3.8,3.3)
\includegraphics[width=3.8cm]{StExa3s1_21.eps}
\put(-3.50,2.70){\scriptsize $\omega = 0$}
\put(-3.50,2.40){\scriptsize $\alpha = \frac{1}{2}\pi$}
\put(-3.50,2.10){\scriptsize $s = 1$}
\end{picture}
\hspace{.5cm}
\begin{picture}(3.8,3.3)
\includegraphics[width=3.8cm]{StExa3s1_22.eps}
\put(-3.50,2.70){\scriptsize $\omega = \frac{1}{4}$}
\put(-3.50,2.40){\scriptsize $\alpha = \frac{1}{2}\pi$}
\put(-3.50,2.10){\scriptsize $s = 1$}
\end{picture}
\hspace{.5cm}
\begin{picture}(3.8,3.3)
\includegraphics[width=3.8cm]{StExa3s1_23.eps}
\put(-3.50,2.70){\scriptsize $\omega = \frac{1}{3}\sqrt{2}$}
\put(-3.50,2.40){\scriptsize $\alpha = \frac{1}{2}\pi$}
\put(-3.50,2.10){\scriptsize $s = 1$}
\end{picture}
\caption{ \small \label{Fig:StExa3s1}
Stability domains for $s=1$. Methods from Example~\ref{Exa3},
with the function $r = r_B$ [eq.~(\ref{eq:rB})] for
$\theta = 1 - \frac{1}{2}\sqrt{2}$ and
$\nu = 2\theta(\hat{a}_{32}-b_2)$, $b_2 = -\frac{1}{2} + \frac{1}{4\theta}$,
$\hat{a}_{32} = \frac{1}{2} + \omega$.
From left to right: $\omega = 0 $ [left],
$\omega = \frac{1}{4}$ [middle] and
$\omega = \frac{1}{3}\sqrt{2}$[right].
Top row for angle $\alpha = 0$, bottom row for $\alpha = \frac{1}{2}\pi$.
For reference: the dash-dotted lines indicate the stability boundaries of the
explicit methods from eq.~(\ref{eq:PPexa1}), and the dashed lines give
the necessary condition $|\phi_B(z_0)| \le 1$, with $\phi_B$ given by 
eq.\ (\ref{eq:NesCondB2}).
[Note: $\phi_B=0$ for $\omega=0$.]
}
\vspace{-3mm}
\end{center}
\end{figure}

\bigskip\noindent
{\it Methods (\ref{eq:Exa3}) with $s=1$\/}:
As a final example we consider the methods from Example~\ref{Exa3} with
$s=1$. With these methods we have $\omega$ as free parameter, and
$\theta = 1 - \frac{1}{2}\sqrt{2}$ is fixed.
Since it concerns here type-{\small B} methods, the stability
domains are empty for $s=2$ or larger. In contrast to this, we do get
rather large domains for $s=1$, as shown in Figure\,\ref{Fig:StExa3s1}.

In this figure, along with the domains ${\cal D}_\alpha$, also the
curves are drawn that are obtained from the necessary condition in
Proposition~\ref{Pro:TypeAlim}(a) for $s=1$. Stability of
the explicit method can be viewed as another necessary condition,
for the special case $z_1=0$. Together these two necessary conditions
give a quite good inclusion for ${\cal D}_\alpha$ for the parameters
$\omega = \frac{1}{4}$ and $\omega = \frac{1}{3}\sqrt{2}$.
For $\omega = 0$ the function $\phi_B$ from (\ref{eq:NesCondB2}) is
equal to zero, so then the necessary condition from
Proposition~\ref{Pro:TypeAlim}(a) vanishes.

\begin{Rem} \rm
For parabolic problems with mixed derivatives on Cartesian grids, one 
can apply dimension splitting with explicit treatment of the mixed 
derivatives. Stability results for such problems can be found in 
\cite{HoMi13,HoWe09}, for example, together with applications for 
option pricing in financial mathematics.
\end{Rem}

\section{Numerical illustration}

In this section we present some numerical test results for a 2D 
reaction-diffusion problem. This problem will be examined with $s=1$ 
and $s=2$ to illustrate the differences between the type-{\small A} 
and type-{\small B} methods.

The test results will be presented for the following methods:
\begin{center}
\begin{tabular}{l}
{SCM-A1} : type-{\small A} method, (\ref{eq:Exa1}) with 
$\theta = 1 - \frac{1}{2}\sqrt{2}$,
\\[1mm]
{SCM-A2} : type-{\small A} method, (\ref{eq:Exa1}) with 
$\theta = \frac{1}{2} + \frac{1}{6}\sqrt{3}$,
\\[1mm]
{SCM-B1} : type-{\small B} method, (\ref{eq:Exa3}) with 
$\theta = 1 - \frac{1}{2}\sqrt{2}$, \,$\omega = 0$,
\\[1mm]
{SCM-B2} : type-{\small B} method, (\ref{eq:Exa3}) with 
$\theta = 1 - \frac{1}{2}\sqrt{2}$, \,$\omega = \frac{1}{3}\sqrt{2}$.
\end{tabular}
\end{center}
The tests were also performed with the type-{\small A} methods 
(\ref{eq:Exa2}), as well as with the variant with $\kappa = \frac{1}{2}$, 
but the results with these methods differed only very slightly from 
the ones with the methods (\ref{eq:Exa1}). In the error plots this would have 
led to lines and markers visually coinciding with those for (\ref{eq:Exa1}).

\subsection{A reaction-diffusion problem with pattern formation}
As test problem we consider the so-called Schnackenberg model
for the interaction of two chemical species, consisting of the 
following system of reaction-diffusion equations
\eq
\setlength{\arraycolsep}{1mm}
\label{eq:Schnak}
\begin{array}{ccl}
u_t &=& D_1(u_{xx}+u_{yy}) \,+\, \kappa(a-u+u^2v) \,,
\\[2mm]
v_t &=& D_2(v_{xx}+v_{yy}) \,+\, \kappa(b-u^2v) \,.
\end{array}
\eeq
This system is considered on the spatial domain $\Omega = [0,1]^2$ and 
time interval $[0,T]$. 
The initial condition is
$$
u(x,y,0) = a + b + 10^{-3}
e^{-100\left(\left(x-\frac{1}{4}\right)^2+\left(y-\frac{1}{6}\right)^2\right)},
\qquad
v(x,y,0) = {b}/{(a+b)^2},
$$
and at the boundaries homogeneous Neumann conditions are imposed.
The parameter values are $D_1=0.05$, $D_2=1$, $\kappa=100$, $a=0.1305$
and $b=0.7695$. The initial condition consists of a small Gaussian
perturbation added to the chemical steady state $u \equiv a+b = 0.90$,
$v \equiv b/(a+b)^2 = 0.95$.
This small perturbation is then amplified and spread, leading to the formation
of patterns with spots.

The space derivatives are approximated by standard
second-order finite differences on a uniform Cartesian mesh.
This problem was used as numerical test in \cite{HuVe03} with
\mbox{$s=1$} and $F_1$ the discretized 2D diffusion operator.
Here we will also consider \mbox{$s=2$} with dimension splitting,
such that $F_1$ contains the discretized $x$-derivatives and $F_2$
the $y$-derivatives.
Computationally the problem becomes much easier with this dimension 
splitting, because all linear systems to be solved are then essentially
three-diagonal.
Tests in \cite{AHHP17} for this problem on a hexagonal
spatial region with domain decomposition splitting revealed that the
errors at output time $T=\frac{1}{2}$ can be quite different from those
at $T=1$.  Therefore the errors will be presented here for both these output
times.

\begin{figure}[t]
\setlength{\unitlength}{1cm}
\begin{center}
\begin{picture}(5,4.2)
\includegraphics[width=5cm]{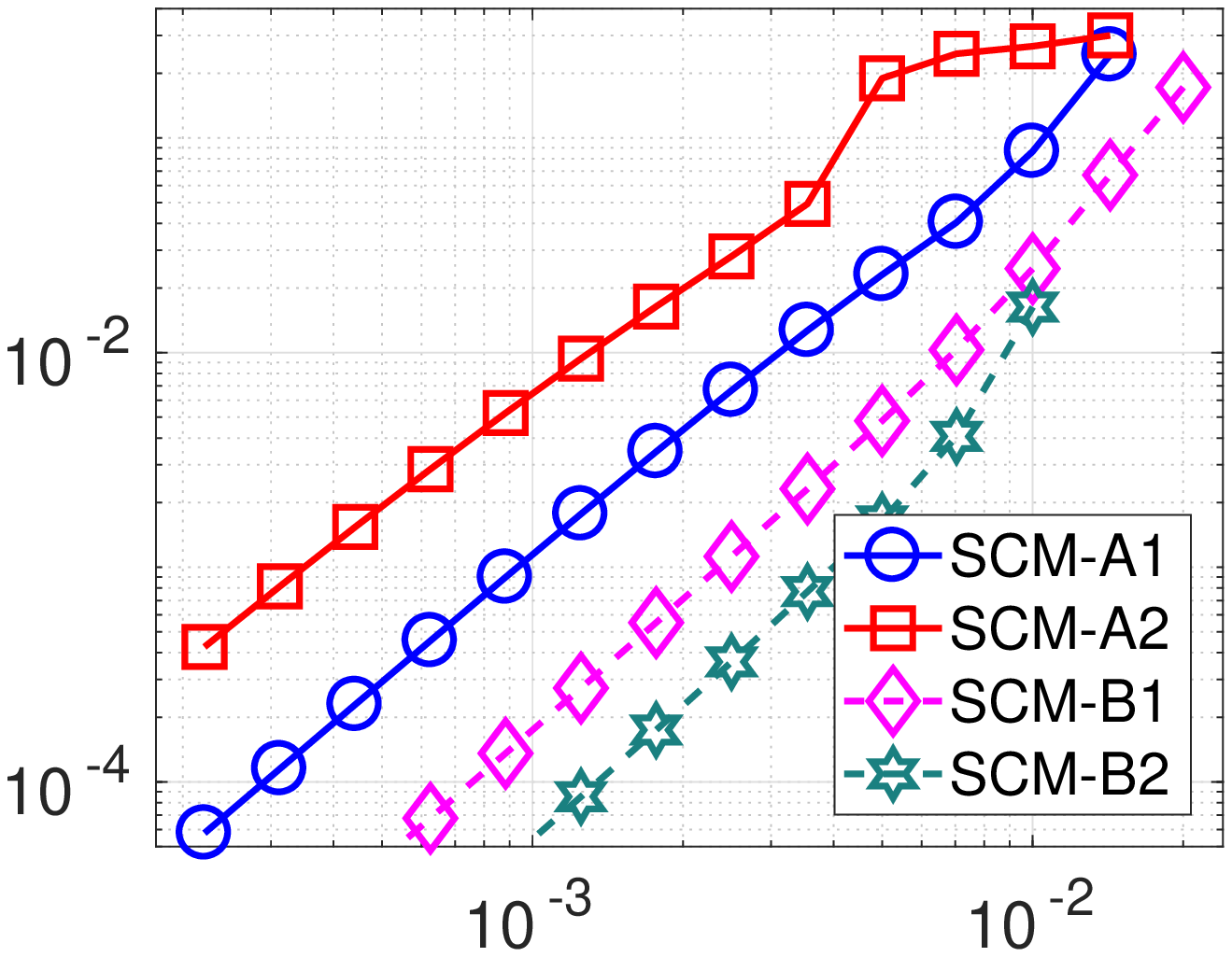}
\put(-4.3,3.5){\scriptsize $L_2$-errors at $T=\frac{1}{2}$}
\put(-4.3,3.2){\scriptsize $s = 1$}
\end{picture}
\hspace{1.5cm}
\begin{picture}(5,4.2)
\includegraphics[width=5cm]{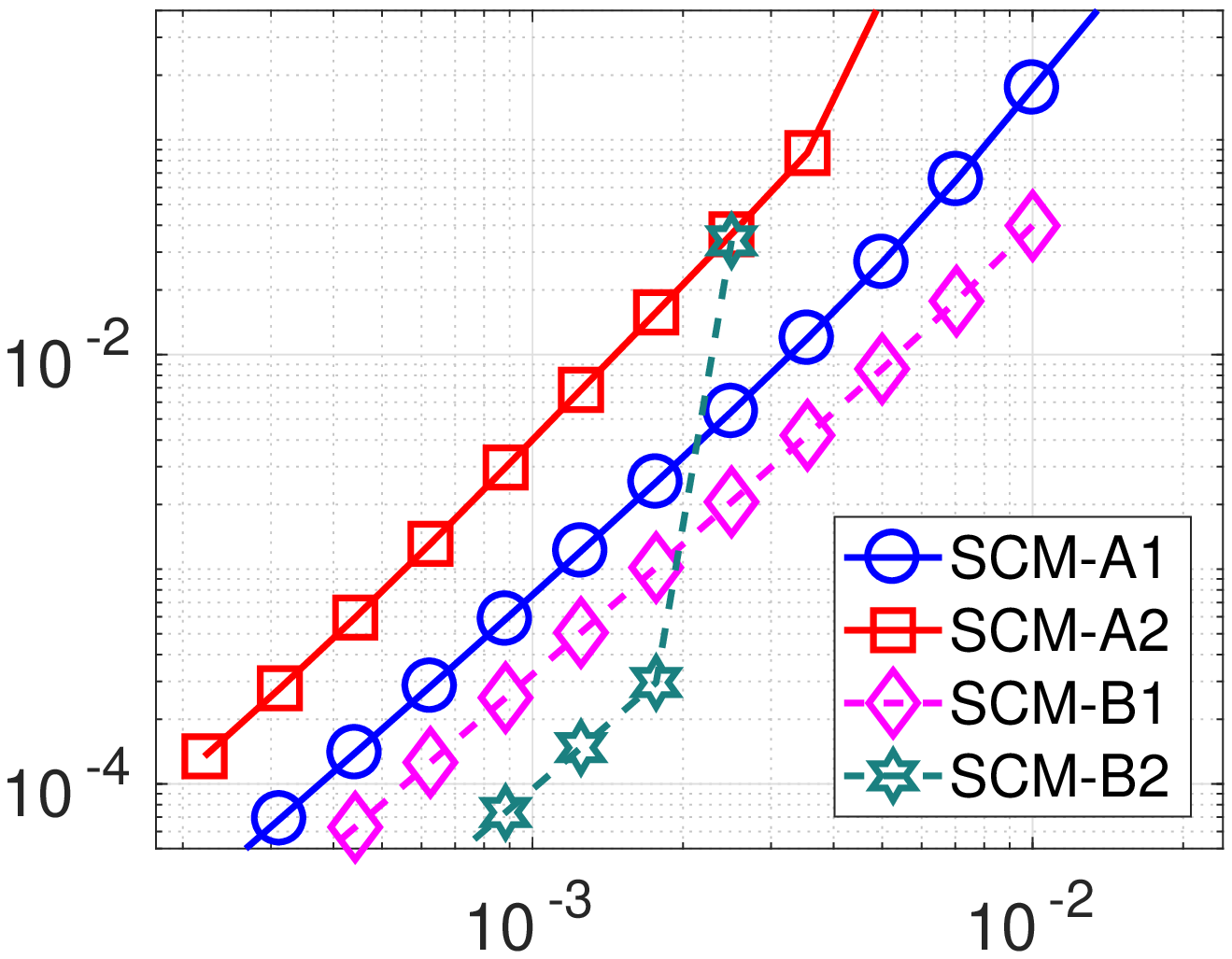}
\put(-4.3,3.5){\scriptsize $L_2$-errors at $T=1$}
\put(-4.3,3.2){\scriptsize $s = 1$}
\end{picture}
\vspace{-3mm}
\caption{ \small  \label{Fig:Resu1a}
Problem (\ref{eq:Schnak}) with $s=1$: $L_2$-errors versus $\dt$ at
output time $T=\frac{1}{2}$ [left] and $T=1$ [right]
with $100\times100$ spatial grid and $\dt = 1/N$,
$N = 50\cdot2^{(j-1)/2}$ rounded to even numbers ($j=1,2,\ldots,14$).
The methods are unstable with $\dt=\frac{1}{50}$, except
for SCM-B1 with $T=\frac{1}{2}$; if $T=1$, method SCM-B2 already appears
unstable with step-size $\dt=\frac{1}{400}$.
}
\end{center}
\vspace{-3mm}
\end{figure}

\begin{figure}[t]
\setlength{\unitlength}{1cm}
\begin{center}
\begin{picture}(5,4.2)
\includegraphics[width=5cm]{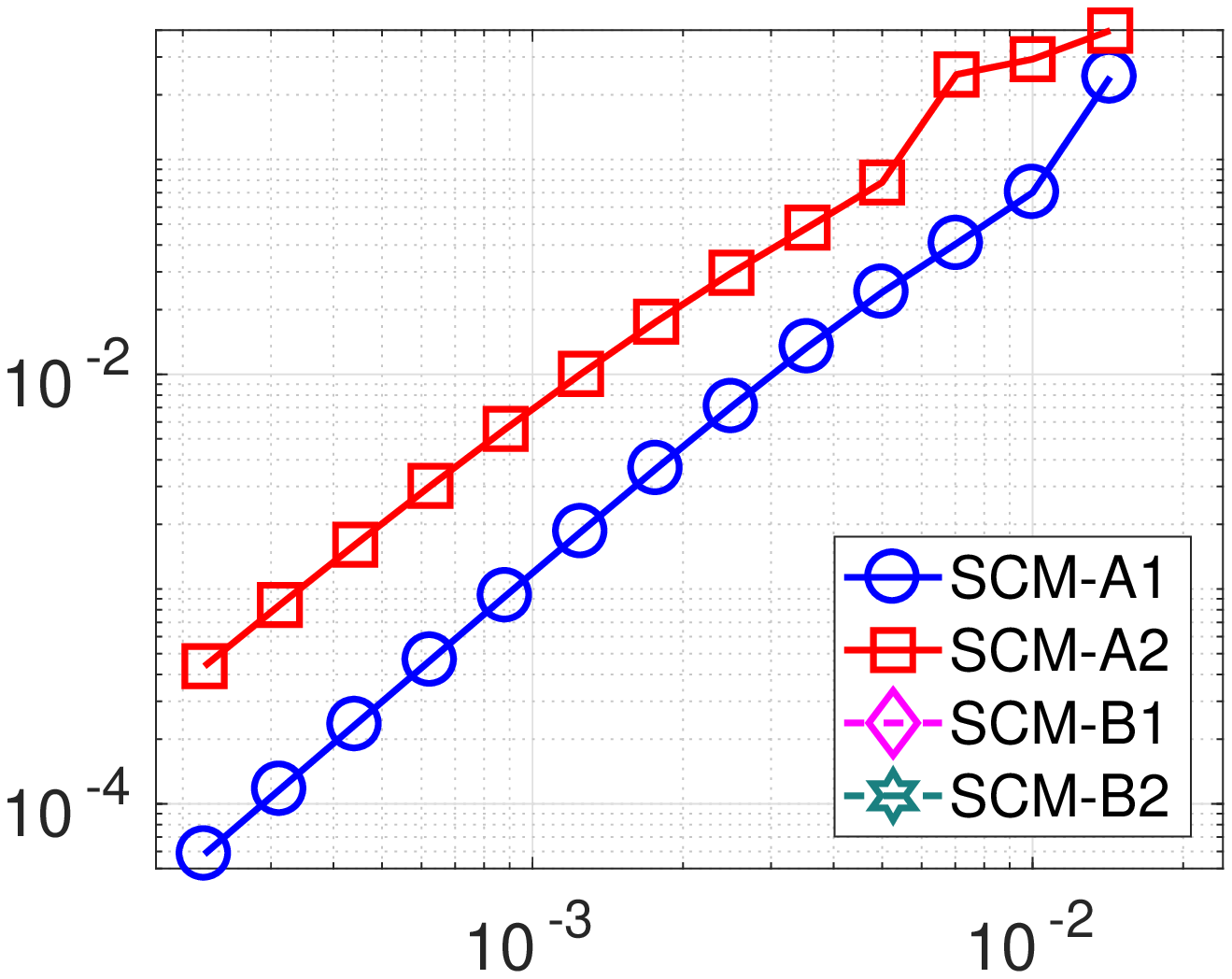}
\put(-4.3,3.5){\scriptsize $L_2$-errors at $T=\frac{1}{2}$}
\put(-4.3,3.2){\scriptsize $s = 2$}
\end{picture}
\hspace{1.5cm}
\begin{picture}(5,4.2)
\includegraphics[width=5cm]{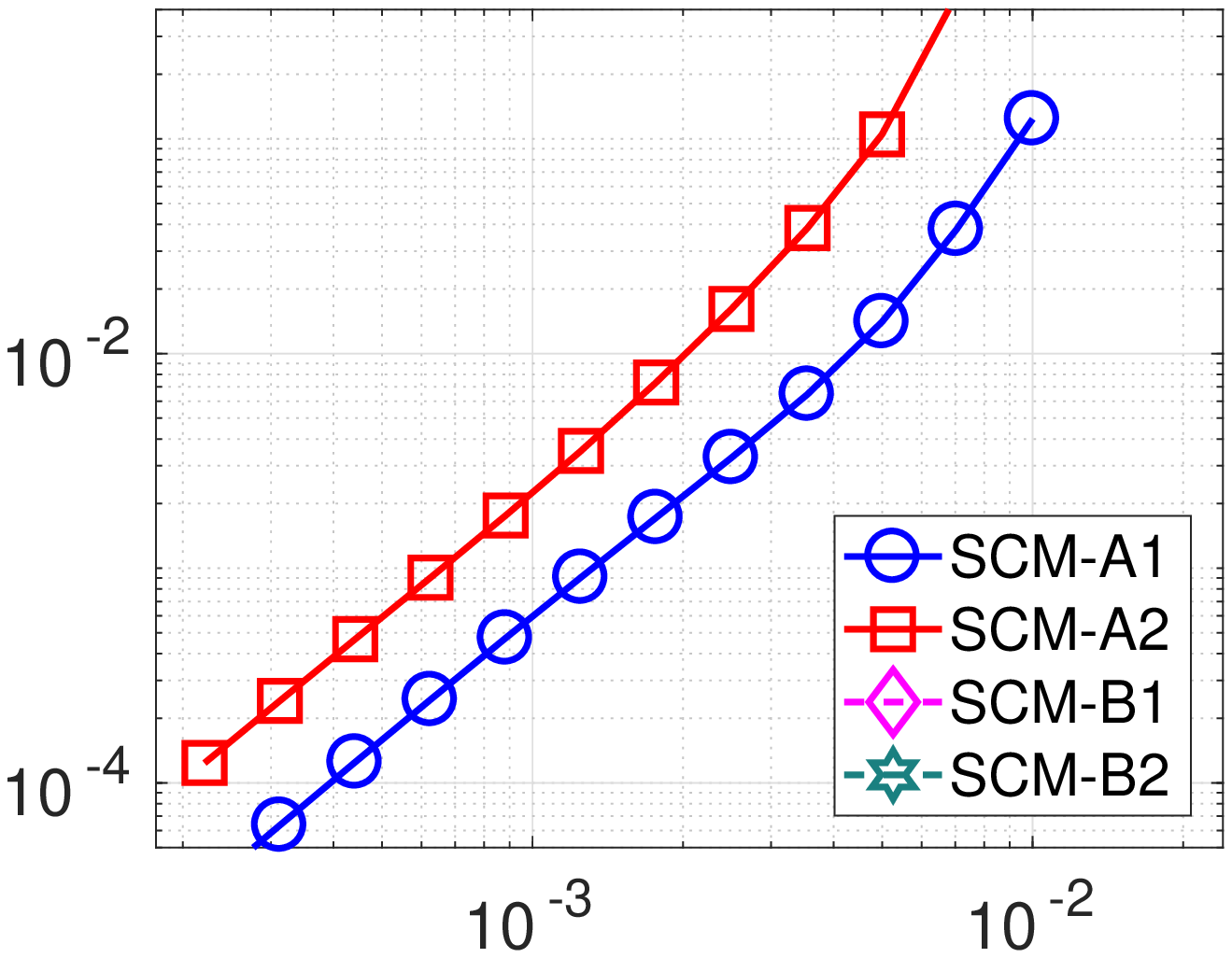}
\put(-4.3,3.5){\scriptsize $L_2$-errors at $T=1$}
\put(-4.3,3.2){\scriptsize $s = 2$}
\end{picture}
\vspace{-3mm}
\caption{ \small  \label{Fig:Resu1b}
Problem (\ref{eq:Schnak}) with $s=2$:
$L_2$-errors versus $\dt$ at output time $T=\frac{1}{2}$ [left] and
$T=1$ [right] with $100\times100$ spatial grid and $\dt = 1/N$,
$N = 50\cdot2^{(j-1)/2}$ rounded to even numbers ($j=1,2,\ldots,14$).
The methods SCM-A1, SCM-A2 are unstable with $\dt=\frac{1}{50}$; the methods
SCM-B1, SCM-B2 are unstable for all the step-sizes.
}
\end{center}
\vspace{-3mm}
\end{figure}

The results with a uniform spatial $100\times100$ grid are shown in
Figure\,\ref{Fig:Resu1a} for $s=1$ and Figure\,\ref{Fig:Resu1b} for $s=2$.
In these figures the temporal $L_2$-errors for the $u$-component are
plotted as function of the step-size.
To compute these temporal errors, a reference solution was found
on this fixed spatial grid by using a very small step-size.
Unstable results are not shown in these error plots.
Consequently, the type-{\small B} methods are listed in the
legends in Figure\,\ref{Fig:Resu1b} even though these methods turned
out to be unstable for $s=2$, in agreement with
Proposition~\ref{Pro:TypeBlim}.

As seen in the figures, the type-{\small A} methods can be used for
the case of dimension splitting with $s=2$. In fact, for this problem
the dimension splitting hardly leads to an increase of the errors.
The type-{\small B} methods, on the other hand, can only be used with
$s=1$, but for that case the errors can be significantly smaller than
those of the type-{\small A} methods.

With respect to stability for $s=1$, the type-{\small B} method with
$\omega=0$ does allow larger steps than the method with
$\omega = \frac{1}{3}\sqrt{2}$, in agreement with the stability
domains shown in Figure\,\ref{Fig:StExa3s1}.
However, once the step-size is small enough to have stability, the errors
for the method with $\omega = \frac{1}{3}\sqrt{2}$ are smaller than
for $\omega=0$.

\section{Concluding remarks}

\subsection{Generalizations}
The general stabilizing correction procedure can be generalized
by treating $F_0$ different from the other $F_j$ in the prediction
steps. Based the formulas (\ref{eq:stageEx}) and (\ref{eq:stageIm})
for the explicit and the implicit method, with coefficients
$\hat{a}_{ik}$ and $a_{ik}$ respectively, we can proceed in the
following way:
\eq
\label{eq:StageSCgen}
\setlength{\arraycolsep}{1mm}
\left\{
\begin{array}{rcl}
x_{i,0} &=&
\displaystyle
u_{n} \,+\, \dt \sum_{k=1}^{i-1} \Big( \hat{a}_{ik} F_0(t_n+c_k\dt, y_k)
\,+\, \sum_{j=1}^{s} \check{a}_{ik} F_j(t_n+c_k\dt, y_k) \Big)
 \,,
\\[2mm]
x_{i,j} &=&
\displaystyle
x_{i,j-1} \,+\, \dt \sum_{k=1}^{i-1} \big(a_{ik} - \check{a}_{ik}\big)
F_j(t_n+c_k\dt, y_k)
\\[6mm]
&  &
\qquad\qquad
\,+\,\, \dt\,a_{ii} F_j(t_n+c_i\dt, x_{i,j})
\qquad\qquad (j = 1,2,\ldots,s) \,,
\\[2mm]
y_i &=& x_{i,s} \,.
\end{array}
\right.
\eeq
Here the new coefficients $\check{a}_{ik}$ should be such that
\eq
\sum_{k=1}^{i-1} \check{a}_{ik} \,=\, \sum_{k=1}^{i-1} \check{a}_{ik}
\eeq
in order to preserve the internal consistency of the scheme, by which 
all vectors $x_{i,j}$, $j=0,1,\ldots,s$, are to be consistent
approximations to $u(t_n+c_i\dt)$.

By allowing such generalizations one can include for example the
modified Craig-Sneyd methods of \cite{HoWe09} and the modified
Douglas method of \cite{AHHP17}.
However, the general formulas that can be obtained this way contain 
new free parameters $\check{a}_{ik}$ for which it is not easy to 
make suitable a priori choices. 

This last point also applies for Runge-Kutta pairs with three or
more stages. It seems, at present, no attempts have been made in
that direction.  Suitable pairs of methods might be found within the
classes of IMEX Runge-Kutta methods derived in \cite{ARS97,BoPa17,KeCa03},
for example.

Other generalization are possible by considering multistep methods.
It is then easy to obtain methods of high order, but stability
becomes more problematic.  Application of such methods to special
classes of parabolic problems from mathematical finance will be
discussed in a separate report.

\subsection{Conclusions}
In this technical note splitting methods have been derived, using the
idea of stabilizing corrections, starting from suitable pairs of 
explicit and diagonally-implicit Runge-Kutta methods with two stages.
In the resulting splitting methods all internal stages provide fully
consistent approximations to the exact solutions. Consequently, steady 
states of the ODE system are preserved as steady states of the splitting 
methods.
  
Linear invariance properties can be ensured in the splitting methods
by performing a finishing stage involving the whole function $F$
(the type-{\small B} methods).
However, it was found that for multiple implicit terms, $s>1$,
unconditional stability properties are then lost.

For the derived splitting methods, linear stability properties have been
studied and some numerical tests were performed, but a full convergence
analysis is still lacking. For practical relevance, such an analysis
should be valid for (semi-discrete) PDEs and stiff ODEs, but this will
make the analysis rather complicated because the local errors must be
studied together with the error propagation, see for example the results
obtained in \cite{Hu02} for $s=1$ and \cite{HoWy15} for $s=2$ with some
{type-{\small A}} methods.

In the numerical example it was seen that type-{\small B} methods may
perform well for $s=1$, but the additional finishing stage makes these
methods unsuitable for problems with multiple implicit terms.
Among the type-{\small A} methods with $\theta$ fixed,
the influence of the second parameter $\kappa$ was marginal in the
tests, producing error plots with almost identical lines.
Stronger nonlinearities might be needed to see significant          
differences.

\appendix

\end{document}